\documentclass[12pt,twoside]{article}

\usepackage{graphicx}
\date{}

\def\sq{{\rm \sqcap \kern-.65em _{-} }}
\centerline{\bf }

\begin{document}

\begin{center}
{\bf  APPROXIMATE IDENTITY AND ARENS REGULARITY OF SOME BANACH ALGEBRAS}
\end{center}

\begin{center}
KAZEM HAGHNEJAD AZAR  AND ABDOLHAMID RIAZI
\end{center}

\small{\noindent{\footnotesize{ABSTRACT. Let $A$ be a Banach algebra with the second dual
$A^{**}$.
 If $A$ has a bounded approximate identity $(=BAI)$, then
$A^{**}$ is unital if and only if $A^{**}$ has a
 $weak^*~bounded~approximate~$$identity(=W^*BAI)$. If $A$ is Arens regular and $A$ \noindent has a
 BAI,
then $A^*$ factors on both sides.
  In this paper we introduce  new concepts $LW^*W$ and $RW^*W$- property
  and we show that under certain conditions if $A$ has $LW^*W$ and $RW^*W$- property, then $A$ is Arens regular
  and also if $A$ is Arens regular,  then
  $A$ has $LW^*W$ and $RW^*W$- property. We also offer some applications of
  these
   new concepts for the special algebras $l^1(G),~L^1(G),~M(G)$, and
   $A(G)$.}}\\\\
\centerline{} {\bf Mathematics Subject Classification:}
 46L06; 46L07; 46L10; 47L25

{\bf  Preliminaries}\\

Through of  this paper, $A$ is  a Banach algebra and $A^*$,
$A^{**}$, respectively, are the first and second dual of $A$. We
say that a bounded net $(e_{\alpha})_{{\alpha}\in I}$ in $A$ is an
approximate identity $(=BAI)$ if,
 for each $a\in A$, $ae_{\alpha}\longrightarrow a$ and  $e_{\alpha}a\longrightarrow a$. For $a\in A$
 and $f\in A^*$ , we denote by $f.a$
 and $ao f$ respectively, the functionals on $A^*$ defined by $<f.a,b>=<f,ab>=f(ab)$ and $<ao f,b>=<f,ba>=f(ba)$.
As since the Banach algebra $A$ can be considered a subspace of its
second dual, so we can define $<a,f>$ as $<f,a>$ for every $a\in
A$ and $f\in
A^*$.\\
 We denote by $A^*A$ and $AA^*$, respectively, $\{f.a:~a\in A~ and ~f\in
 A^*\}$ and
  $\{ao f:~a\in A ~and ~f\in A^*\}$ clearly these two sets are subsets of $A^*$.
 We say that the Banach algebra $A$ is unital
if there exists an element   $e\in A$ such that $e.x=x.e=x$ for
each $x \in A$.\\
 Let $A$ be a Banach algebra with a $BAI$. If the
equality $A^*A=A^*,~~(AA^*=A^*)$ holds, then we say that $A^*$
factors on the left (right). If both equalities $A^*A=AA^*=A^*$
hold, then we say
that $A^*$  factors on both sides.\\
Arens [1] has shown that given any Banach algebra $A$, there exist
two algebra multiplications on the second dual  of $A$ which
extend multiplication on  $A$. In the following, we introduce both
multiplication which are given in [11].\newpage
 Let $a,b \in A$
and $f\in A^*$  and  $F,G \in A^{**}$ then the first Arens
multiplication is  defined as
$$<f.a,b>=<f,ab>,$$
$$<F.f,a>=<F,f.a>,$$
$$<F.G,f>=<F,G.f>.$$
Clearly $F.f\in A^*$ and $F.G\in A^{**}$. We use the notions
$(A^{**},.)$ for $A^{**}$ equipped  with first Arens
multiplication also we use $fa$ instead of $f.a$ for all $a\in A$
and $f\in A^*$.
 \\
The second Arens product is defined as follows \\
For $a,b\in A$ , $f\in A^*$ and $F,G\in A^{**}$  the elements $ao
f$ , $fo F$ of $A^*$ and $Fo G$ of $A^{**}$ are defined
respectively  by the equalities
$$<ao f,b>=<f,ba>,$$
$$<fo F,a>=<F,ao f>,$$
$$<Fo G,f>=<G,fo F>.$$
An element $E$ of $A^{**}$ is said to be a mixed unit if $E$ is a
right unit for the first Arens multiplication and a left unit for
the second Arens multiplication. That is, $E$ is a mixed unit if
and only if,
for each $F\in A^{**}$, $F.E=Eo F=F$.\\
We say that $A^{**}$ is unital with respect to the first Arens
product, if there exists an element $E\in A^{**}$ such that
$F.E=E.F=F$ for all $F\in A^{**}$ , and $A^{**}$ is unital with
respect to the second Arens product , if there exists an element
$E\in A^{**}$ such that $Fo E=Eo F=F$ for all $F\in A^{**}$. By
[3, p.146], an element $E$ of $A^{**}$  is  mixed
      unit if and only if it is a $weak^*$ cluster point of some BAI $(e_\alpha)_{\alpha \in I}$  in
      $A$.\\
We say that $A^{**}$ has a
$weak^*~bounded~left~approximate~identity(=W^*BLAI)$ with respect
to the first Arens product, if there is a bounded net as
$(e_{\alpha})_{\alpha}\subseteq A$ such that for all $F\in
A^{**}$ and $f\in A^*$, we have
 $<e_{\alpha}.F,f>\rightarrow<F,f>$. The definition of $W^*RBAI$ is similar to $W^*LBAI$   and
 if $A^{**}$ has both
$W^*LBAI$ and  $W^*RBAI$, then we say that $A^{**}$ has $W^*BAI$.\\
 Suppose that $F,G\in A^{**}$ and $F.G,~Fo G$ are the first and second Arens Multiplications in $A^{**}$, respectively.
 Then the mapping $F \rightarrow F.G $, for $G$  fixed in $ A^{**}$, is $weak^*-weak^*$
 continuous, but the mapping $F \rightarrow G.F $ for $G$ fixed in  $ A^{**}$  is not in general
  $weak^*-weak^*$ continuous on $A^{**}$ unless $G\in A$.\\
  As an example for the algebras $L^1(G)^{**}$ and $M(G)^{**}$
  whenever
  $G$ is an infinite topological group the mapping $F\rightarrow G.F$, in
  general, is not $~weak^*-weak^*~
   continuous~on$ $L^1(G)^{**}$ and $M(G)^{**}$, see [10,11].
  Hence, the first topological center of $ A^{**}$
  with respect to first Arens  product is defined as follows
   $$Z_1=\{G \in A^{**}:  ~ F \longrightarrow G.F ~is ~w^*-w^*-
   continuous~on
   ~A^{**}\}.$$

Also, for the second Arens multiplication in $A^{**}$ we know that
$ao f , fo F\in A^*$ and $Fo G\in A^{**}$. For $G$ fixed in
$A^{**}$ , the mapping $F\rightarrow Go F$ is $weak^*-weak^*$
continuous on $A^{**}$, but the mapping $F\rightarrow Fo G$ is not
in general $weak^*-weak^*$ continuous on $A^{**}$ unless $G\in A$.
  Whence, the second topological center of $ A^{**}$
  with respect to second Arens product is defined as follows
 $$Z_2=\{G \in A^{**}: ~ F \longrightarrow Fo G ~ is~
w^*-w^* -continuous ~on~  A^{**}\}.$$
 It is clear that $A\subseteq
Z_1\bigcap Z_2$ and $Z_1,~Z_2$ are closed subalgebras of $A^{**}$
endowed with the first
 second Arens multiplication, respectively.

If, for each $F,G \in A^{**}$, the equality  $F.G=FoG$  holds,
then the algebra $ A$ is said to be Arens regular, see [1,2]. In this case $Z_1=Z_2=A^{**}$.\\
The other extreme situation is that $Z_1=A$ , in this case $A$ is
called left strongly Arens irregular, see [9,10,14].
\\
We recall that the topological center of $A^{**}$ is defined the
set of all functionals $F\in A^{**}$ which satisfy $F.G=Fo G$  for
all $G\in A^{**}$, see [11]. In other words, the topological
centers of $A^{**}$ with respect to the first and second
Arens products  can also be defined as following sets respectively\\
$$Z_1=\{F\in A^{**} :~F.G=Fo G~~~  \forall G\in A^{**}\},$$

$$Z_2=\{F\in A^{**}:~G.F=Go F~~~  \forall G\in A^{**}\}.$$
For Banach algebra $A$ the topological center of the algebra
$(A^*A)^*$ is defined to be following set, see [11].
$$\widetilde{Z}=\{ \mu\in (A^*A)^*:~\lambda\rightarrow~\lambda.\mu ~~is~
w^*-w^*-continuous~on~(A^*A)^*\}.$$
 For all $a\in A$ , $f\in A^*$ and $F\in A^{**}$
 we have the following statements \\
i) $f.a=fo a$ ~and~ $a.f=ao f$.\\
ii) $a.F=ao F$ ~and~ $F.a=Fo a$.\\
Because if $b\in A$, then we have
$$<fo a,b>=<b, foa>=<ab,f>=f(ab).$$
Also, we have $<f.a,b>=f(ab)$, and we conclude that $f.a=fo a$.\\
Consequently we can write the following relations
$$<a.F,f>=<a,F.f>=<F.f,a>=<F,f.a>=<F,fo a>=<ao F,f>.$$
Thus, we conclude that $a.F=ao F$. \\
In above, part (ii)  shows that $A^{**}$ has $W^*LBAI$ with
respect to the first Arens product if
 and only if it has $W^*LBAI$ with respect to the second Arens product.\\
 A functional $f$ in $A^*$ is said to be $wap$ (weakly almost
 periodic) on $A$ if the mapping $a\rightarrow f.a$ from $A$ into
 $A^{*}$ is weakly compact. Pym in [15] showed that  this definition to the equivalent following condition\\
 For any two net $(a_{\alpha})_{\alpha}$ and $(b_{\beta})_{\beta}$
 in $\{a\in A:~\parallel a\parallel\leq 1\}$, we have\\
$$\lim_{\alpha}\lim_{\beta}<f,a_{\alpha}b_{\beta}>=\lim_{\beta}\lim_{\alpha}<f,a_{\alpha}b_{\beta}>,$$
whenever both iterated limits exist. The collection of all $wap$
functionals on $A$ is denoted by $wap(A)$. Also we have
$wap(A)=A^*$ if and only if $<F.G,f>=<FoG,f>$ for every $F,G \in
A^{**}$, see [19]. \\
In this paper, the notations $WSC$ is used for weakly sequentially
complete Banach algebra $A$ and $WCC$, for weakly completely
continuous Banach algebra $A$, that is, $A$ is said to be weakly
sequentially complete, if every weakly Cauchy sequence in $A$ has
a weak limit, and $A$ is said to be $WCC$, if for each $a\in A$,
the
multiplication operator $x\rightarrow ax$ is weakly compact.\\
In this paper when $A^{**}$ is referred to as an algebra it will
be with respect to the first Arens product unless
the second is mentioned explicitly.\\\\

\textbf{\large{ Main results }}\\\\
{\it{\bf Theorem 1.}} Let $A$ has a BAI $(e_{\alpha})_{\alpha}$. Then the following assertions  hold:\\
i)  For all $f\in A^*$, we have $f.e_\alpha\stackrel{w^*} {\rightarrow} f$ and $e_\alpha of\stackrel{w^*} {\rightarrow} f$.\\
ii)  $A^*$ factors on the right if and only if $A^{**}$ has a $W^*LBAI$.\\
iii)  $A^*$ factors on the left  if and only if $A^{**}$ has a $W^*RBAI$.\\
iv)  If $A^*$ factors on the left and $E\in A^{**}$ such that
$e_{\alpha} \stackrel{w^*} {\rightarrow} E$, then $E$ is unit
element of $A^{**}$. Also, if $A^*$ factors on the right and $E\in
A^{**}$ such that $e_{\alpha} \stackrel{w^*} {\rightarrow} E$,
then $E$ is unit
element of $(A^{**},o)$.  \\
v) If $A^*$ factors on the left,  then $w^*-closure(AA^{**})=A^{**}$.\\\\
 Proof. i) The proof of (i) is clear. \\
  ii) Let $A^*A=A^*$, and $F\in A^{**}$ , $f\in A^*$.
Then, by [11, 2.1] , we have $f.e_{\alpha} \stackrel{w} {\rightarrow} f$, and so\\
$$<e_{\alpha}.F,f>=<e_{\alpha},F.f>=<F.f,e_{\alpha}>=<F,f.e_{\alpha}>\longrightarrow <F,f>.$$
Consequently,  $e_{\alpha}.F \stackrel{w^*} {\rightarrow} F$.\\
Conversely, let $e_{\alpha}.F \stackrel{w^*} {\rightarrow} F$ for all $F\in A^{**}$. Then, for every $f\in A^*$, we have $<e_{\alpha}.F,f>\longrightarrow<F,f>$. Since $<e_{\alpha}.F,f>=<F,f.e_{\alpha}>,$ we have\\
 $<F,f.e_{\alpha}>\longrightarrow <F,f>.$ Therefore by [11,  2.1], we are done.\\
 iii) The proof of (iii) is the same as (ii). \\
iv) Let $e_{\alpha} \stackrel{w^*} {\rightarrow} E$.
 Then, for all $F\in A^{**}$, by using  [11, 2.1], we have the following relations:
$$<E.F,f>=<E,F.f>=\lim_{\alpha}<e_{\alpha},F.f>=\lim_{\alpha}<F,f.e_{\alpha}>\rightarrow<F,f>.$$
Hence, we conclude that $E.F=F$.\\
Now let $f\in A^*$ and $a\in A$. Then , we have
$$<E.f,a>=<E,f.a>=\lim_{\alpha}<e_{\alpha},f.a>=\lim_{\alpha}<f,a.e_{\alpha}>=<f,a>$$
Consequently, we have $<F.E,f>=<F,E.f>=<F,f>$ hence $F.E=F$  which implies that $E$ is a unit element of $A^{**}$.\\
The next assertion is similar.\\
v) Let $A^*A=A^*$, and let $f\in A^*$. By [11, 2.1], we have
$f.e_{\alpha} \stackrel{w} {\rightarrow} f$. For each $F\in
A^{**}$, since $e_{\alpha}.F\in AA^{**}$, we have
$$<e_{\alpha}.F,f>=<F,f.e_{\alpha}>\rightarrow <F,f>.$$
Thus, we conclude that $e_{\alpha}.F \stackrel{w^*} {\rightarrow}
F$
and so $F\in w^*-closure(AA^{**})$. \\\\
By the above theorem and by [11, 2.2.a], for the Banach algebra
$A$ which has $BAI$, we conclude that $A^{**}$ has $W^*LBAI$ if
and only if \\$\{F\in A^{**}:~A.F\subseteq A\}\subseteq Z_1$, and
also $A^{**}$ has $W^*RBAI$ if and only if $\{F\in A^{**}:~Fo
A\subseteq A\}\subseteq Z_2$. For $A=L^1(G)$ where $G$ is locally
compact finite group, we have $\{F\in A^{**}:~A.F\subseteq
A\}\subseteq Z_1=L^1(G)^{**}$ hence $L^1(G)^{**}$ has $W^*LBAI$.
But when $G$ in locally compact infinite group
this relation is not satisfied for $L^1(G)^{**}$ consequently $L^1(G)^{**}$ has not
$W^*LBAI$.\\\\\\
{\it{\bf  Corollary 2.}} Let $A$ has BAI. Then the following statements hold.\\
i) The algebra $(A^{**},.)$ is unital if and only if $A^{**}$ has $W^*LBAI$.\\
ii) The algebra $(A^{**},o)$ is unital if and only if $A^{**}$ has $W^*RBAI$.\\
iii) $A^{**}$ and $(A^{**},o)$ are unital if and only if $A^{**}$ has $W^*BAI$.\\
iv) Suppose that $A$ is a $WSC$ with sequentially $BAI$.
 Then $A^{**}$ has $W^*LBAI$ if and only if $A^{**}$ has
$W^*RBAI$, and also $(A^{**},.)$ is unital if and only if
$(A^{**},o)$ is
unital.\\\\
Proof:  Using  Theorem 1 and [11, 2.2], the proofs of (i) , (ii)
are clear. The proof of (iv), follows from [11, 2.6] and Theorem
1.
  \\\\
As another  application of Theorem 1, let $\Omega$ be  spectrum of
$L^{\infty}([0,1])$.
 Since  $L^{\infty}([0,1])$
 dose not factors on the left and right, then by Theorem 1, $L^1([0,1])^{**}$ has not $W^*LBAI$
 and
  $W^*RBAI$, and so $M(\Omega)=L^1([0,1])^{**}$ has not $W^*LBAI$
  and
  $W^*RBAI$. Consequently by Corollary 2,
   $M(\Omega)$  is not unital  which it follows  that $\Omega$ is not a group. For more
   details
   see
    [9,16 ].
   \\\\
{\it{\bf Corollary 3.}} Let $A$ be  Arens regular with a $BAI$. Assume
also that $A$ is a two sided ideal in $A^{**}$,
then $A^{**}$ has a $W^*BAI$ and so is unital.\\\\
Proof. Use  Theorem 1 and [17, 3.2].\\\\
{\it{\bf Theorem 4.}} Suppose $f\in wap(A)$ and $G\in A^{**}$ such that
$a_{\alpha} \stackrel{w^*} {\rightarrow} G$ where
$(a_{\alpha})_{\alpha}\in A$. Then $f.a_{\alpha} \stackrel{w}
{\rightarrow} f.G$.\\\\
Proof. First we show that $f.a_{\alpha} \stackrel{w^*}
{\rightarrow} f.G$. Let $b\in A$, so we have
$$<f.a_{\alpha} , b>=<a_{\alpha} , b.f>\rightarrow <G,b.f>=<f.G,b>.$$
Thus we have $f.a_{\alpha} \stackrel{w^*} {\rightarrow} f.G$. \\
Since $f\in wap(A)$,~ $<G.F,f>=<GoF,f>$  for all $F,G\in A^{**}$. We have \\
$$<f.G,F>=<G.F,f>=<GoF,f>=<F,foG>=<foG,F>,$$
it follows that $f.G=foG\in A^*$. For every $F\in A^{**}$, there
is $(b_{\beta})_{\beta}\subseteq A$ such that
$b_{\beta}\stackrel{w^*} {\rightarrow}F$, so we have
$$\lim_{\alpha}<F,f.a_{\alpha}>=\lim_{\alpha}\lim_{\beta}<b_{\beta},f.a_{\alpha},>=
\lim_{\alpha}\lim_{\beta}<f,a_{\alpha}b_{\beta}>,$$
$$=\lim_{\beta}\lim_{\alpha}<f,a_{\alpha}b_{\beta}>=\lim_{\beta}<f.G,b_{\beta}>=\lim_{\beta}<b_{\beta},f.G>,$$
$$=<F,f.G>.$$
Thus, we conclude that $f.a_{\alpha} \stackrel{w} {\rightarrow}
f.G$.\\\\
{\it{\bf Corollary 5.}} Suppose $A$ has a BAI $(e_{\alpha})_{\alpha}$
and $f\in wap(A)$. Then $f.e_{\alpha}
\stackrel{w} {\rightarrow} f$.\\\\
Proof. By [5,p.146], there is a mixed unit  $E\in A^{**}$ such
that $e_{\alpha} \stackrel{w^*} {\rightarrow} E$. Since $f\in
wap(A)$, $<E.F,f>=<EoF,f>=<F,f>$ for all $F\in A^{**}$. So, by
Theorem 4, we have $f.e_{\alpha} \stackrel{w} {\rightarrow}
f.E=f$.\\\\
{\it{\bf Corollary 6.}} Suppose $A$ is Arens regular and $A$ has a BAI.
Then we have the following assertions\\
i) $A^*$ factors on both sides.\\
ii) $A^{**}$ is unital, and if $A$ is $WSC$ then $A$ is also
unital.\\\\
Proof. i) Since $A$ is Arens regular, $wap(A)=A^*$ hence by
Corollary 5, we have $f.e_{\alpha} \stackrel{w} {\rightarrow} f$
for all $f\in A^*$. Then, it follows that $A^*$ factors on the
left by [11, 2.1]. Since $A$ is Arens regular, by [11, 2.10] $A^*$
factors on the right and so $A^*$
factors on both sides.\\
 ii) By [11, 2.2], since $A^*$ factors on both
sides, $(A^{**},.)$ and $(A^{**},o)$ are unital. The second claim hold by using [2,  2.6] and part (i).\\\\
 Since $L^1(G)^{*}$, $M(G)^{*}$ and
$A(G)^{*}$ dose not  factors on the left and right whenever $G$ is
infinite group, by Corollary 6 we obtain that $L^1(G)^{**}$,
$M(G)^{**}$ and
$A(G)^{**}$ are not  Arens regular.\\\\
{\it{\bf Theorem 7.}} For Banach algebra $A$ we have the following
assertions\\
i) Assume that $f\in A^*$ and $T_f$ is the linear operator from
$A$ into $A^*$ defined by $T_fa=f.a$. Then, $f\in wap(A)$ if and
only if the adjoint of $T_f$ is
$weak^*-weak$ continuous.\\
ii) If $A$ has a BAI and $\widetilde{Z}=(A^*A)^*$, then
$A^*A\subseteq wap(A)$.\\\\
 Proof. i)  Let $f\in wap(A)$ and $T_f^*:A^{**}\rightarrow A^*$ be
the adjoint of $T_f$. Thus for every $F\in A^{**}$ and $a\in A$,
we
have $<T_f^*F,a>=<F,T_fa>$.\\
Suppose $(F_{\alpha})_{\alpha}\subseteq A^{**}$  such that
$F_{\alpha} \stackrel{w^*} {\rightarrow} F$. We show that
$T_f^*F_{\alpha} \stackrel{w} {\rightarrow} T_f^*F$. Let $G\in
A^{**}$ and $(a_{\beta})_{\beta}\subseteq A$ such that $a_{\beta}
\stackrel{w^*} {\rightarrow} G$. Since $f\in wap(A)$, we have
$<G.F_{\alpha},f>\rightarrow <G.F,f>.$ Hence for fixed $\alpha$ we
have the following relations
$$<G,T_f^*F_{\alpha}>=lim_{\beta}<T_f^*F_{\alpha},a_{\beta}>=lim_{\beta}<F_{\alpha},Ta_{\beta}>
=lim_{\beta}<F_{\alpha},f.a_{\beta}>$$
$$=lim_{\beta}<F_{\alpha}.f,a_{\beta}>=<G,F_{\alpha}.f>,$$
hence
$$lim_{\alpha}<G,T_f^*F_{\alpha}>=lim_{\alpha}<G,F_{\alpha}.f>=lim_{\alpha}<G.F_{\alpha},f>
=lim_{\alpha}<G.F,f>$$
$$=lim_{\alpha}<G,F.f>=<G,T_f^*F>.$$
We conclude that $T_f^*F_{\alpha} \stackrel{w} {\rightarrow}
T_f^*F$, thus $T^*_f$ is $weak^*-weak$ continuous.\\
Conversely, let $T^*$ be $weak^*-weak$ continuous and let
$(F_{\alpha})_{\alpha}\subseteq A^{**}$ be such that $F_{\alpha}
\stackrel{w^*} {\rightarrow} F$. Then we have
$$<G.F_{\alpha},f>=<G,T_f^*F_{\alpha}>\rightarrow
<G,T_f^*F>=<G.F,f>.$$ It follow that $f\in wap(A)$.\\
ii) Let $\widetilde{Z}=(A^*A)^*$ and $F\in A^{**}$. We define
$\mu\in (A^*A)^*$ to be the restriction of $F$ to $(A^*A)^*$.
Since $\widetilde{Z}=(A^*A)^*$, $\mu\in \widetilde{Z}$. Also for
all $a\in A$, we have $a\mu=aF$.\\
By [11, 3.2], since $A\widetilde{Z}=AZ_1$ and $AZ_1\subseteq Z_1$,
we have $aF\in Z_1$. Therefore, for all $G\in A^{**}$ and $f\in
A^*$ we have the following relations
$$<F.G,f.a>=<aF.G,f>=<aFoG,f>=<aFoG,f>=<FoG,f.a>.$$
Consequently $f.a\in wap(A)$.\\\\
 Theorem 7 shows that a Banach algebra  $A$ is Arens
regular if and only if $T^*_f$
is $weak^*-weak$ continuous for every $f\in A^*$.\\\\
 Now   we define new concepts namely
$LW^*W$-property and $RW^*W$-property for Banach algebra $A$ and
compare
them with Arens regularity of Banach algebra $A$.\\\\
{\it {\bf Definition 8.}} Let $(f_{\alpha})_{{\alpha}\in
I}\subseteq A^*$. If for $a\in A$, the convergence
 $a.f_{\alpha} \stackrel{w^*} {\rightarrow} 0$ implies $a.f_{\alpha} \stackrel{w} {\rightarrow} 0$,
 then we say that $a$ has $Left-Weak^*-Weak$
  property or $LW^*W$-property with respect to the first Arens product. The definition of $RW^*W$-property
  is similar. \\
We say that $A$ has the $LW^*W$-property if, for every $a\in A$,
$a$ has $LW^*W$-property.
\\\\
{\it{\bf   Remark and Example 9.}} We know that if $A$ is reflexive,
then $weak$ and $weak^*$ topologies coincide consequently $A$ has
$LW^*W$ and $RW^*W$-property, so it is clear that if $A$ has a
member which has not $LW^*W$-property or $RW^*W$-property, then
$A$ is not reflexive. If  $A$ is commutative, then it is clear
that the $LW^*W$-property and $RW^*W$-property are equivalent.
For Banach algebras $l^1(G)$ and $M(G)$ where $G$ is an infinite
group, the unit elements of them have not $LW^*W$-property or
$RW^*W$-property
(We know that ${l}^{\infty }(G)^*\neq l^1(G)$ and  $M(G)^{**}\neq M(G)$).\\
For a Banach algebra $A$, if  $A^{**}=A^{**}c$ or $A^*=cA^*$  for
some $c\in A$ and $c$ has $LW^*W$-property, then it is easy to see
that $ A=A^{**}$, and also if  $A^{**}=cA^{**}$ or $A^*=A^*c$ for
some $c\in A$ and $c$ has $RW^*W$-property, then $ A=A^{**}$.\\
Now we give some examples of the Banach algebras that have
$LW^*W$-property or $RW^*W$-property.\\
 Let $S$ be a compact semigroup and $A=C(S)$. Then every $f\in
C(S)$ such that $f\geq \alpha $ where $\alpha\in (0,\infty)$ has
$LW^*W$-property. For prove this, let
$(\mu_{\alpha})_{\alpha}\subseteq C(S)^*=M(S)$ and
$\mu_{\alpha}.f\stackrel{w^*} {\rightarrow} 0$ hence
$<\mu_{\alpha}.f,g>\rightarrow 0$ for all $g\in C(S)$.
Consequently we have
$$<\mu_{\alpha}.f,g>=<\mu_{\alpha},f.g>=\int_Sfgd\mu_{\alpha}\rightarrow
0.$$ If we set $g=\frac{1}{f}$, then $\mu_{\alpha}\rightarrow 0$.
Now let $F\in C(S)^{**}=M(S)^*$. Then we have
$$<F,\mu_{\alpha}.f>=<F.\mu_{\alpha},f>=\int_SfdF\mu_{\alpha}\rightarrow
0.$$ It follows that $\mu_{\alpha}.f\stackrel{w} {\rightarrow}
0$.\\
 In Corollary 15, we will show that if the Banach algebra $A$
is Arens regular, then $A$ has  both $LW^*W$-property and
$RW^*W$-property.\\\\
In the following, we design some condition for a Banach algebra
$A$ with $LW^*W$-property
or $RW^*W$-property, so such Banach algebras are Arens regular.\\\\
{\it{\bf Theorem 10.}} For Banach algebra $A$ the following statements  hold.\\
i)   If $A^{**}=cA^{**}$ for some $c\in A$ and $c$ has $LW^*W$-property, then $A$ is Arens regular.\\
ii) If $A^{**}=A^{**}c$ for some $c\in A$ and $c$ has
$RW^*W$-property, then
$A$ is Arens regular.\\
iii) If $A^*=A^*c$ and $c$ has $LW^*W$-property, then
$A^{**}c\subseteq Z_1\cap Z_2$ and $cA^*\subseteq wap(A)\subseteq A^*c$.\\
iv) If $A^*=cA^*$ and $c$ has $RW^*W$-property, then
$cA^{**}\subseteq Z_1\cap Z_2$ and $A^*c\subseteq wap(A)\subseteq cA^*$.\\\\
Proof.  i)
  Let $F,G\in A^{**}$.  We show
that the mapping $F\rightarrow G.F$ is $w^*-w^*-continuous$. Let
$f\in A^*$ and $(F_{\alpha})_{\alpha}\subseteq A^{**}$ be such
that $F_{\alpha} \stackrel{w^*} {\rightarrow} F$. Then for all
$a\in A$ we have
$$<F_{\alpha}.f,a>=<F_{\alpha},f.a>\rightarrow<F,f.a>=<F.f,a>.$$
Consequently, we have $F_{\alpha}.f \stackrel{w^*} {\rightarrow}
F.f$ for all $f\in A^*$ which implies that
 $(F_{\alpha}.f )a \stackrel{w^*} {\rightarrow} (F.f)a$ for all $a\in A$.
 Since $c$ has $LW^*W$-property, $(F_{\alpha}.f )c \stackrel{w} {\rightarrow} (F.f)c$.
 Thus for all $G\in A^{**}$, we have
$$<cG,F_{\alpha}.f>=<G,(F_{\alpha}.f)c>\rightarrow<G,(F.f)c>=<cG,F.f>.$$
Since $A^{**}=cA^{**}$, we can replace $cG$ by $G$, and
consequently we have \\$<G.F_{\alpha},f>\rightarrow<G.F,f>$ for
all $f\in A^*$. Therefore $G\in Z_1$  so that $ Z_1=A^{**}$.\\
ii) The proof of (ii) is similar to (i).\\
iii) We show that the mapping $F\rightarrow Gc.F$ is
$w^*-w^*$-continuous whenever $F,G\in A^{**}$. Let
$(F_{\alpha})\subseteq A^{**}$ be a net such that $F_{\alpha}
\stackrel{w^*} {\rightarrow} F$. Then, we have $(F_{\alpha}.f )
\stackrel{w^*} {\rightarrow} (F.f)$ for every  $f\in A^*$, hence
$(F_{\alpha}.f-F.f)  \stackrel{w^*} {\rightarrow}0 $. Since
$(F_{\alpha}.f-F.f)\in A^*=A^*c$, there is
$(g_{\alpha})_{\alpha}\subseteq A^*$ such that
$(F_{\alpha}.f-F.f)=g_{\alpha}c$. Then we have
$c(g_{\alpha}c)\stackrel{w^*} {\rightarrow}0$, it follow that
$c(g_{\alpha}c)\stackrel{w} {\rightarrow}0$, since $c$ has
$LW^*W$-property. Consequently for all $G\in A^{**}$, we have
$<G,c(g_{\alpha}c)>\stackrel{} {\rightarrow}0$ which implies
$<Gc,(F_{\alpha}.f-F.f)>\rightarrow 0$ hence $Gc\in Z_1$ so
$A^{**}c\subseteq Z_1$. Similarly $A^{**}c\subseteq Z_2$. Hence we
have $A^{**}c\subseteq Z_1\cap
Z_2$.\\
Now let $F,G\in A^{**}$. Since $Gc\in Z_2$, we have
$Fo(Gc)=F.(Gc)$, consequently for all $f\in A^*$, we have the
following statements
$$<FoG,cf>=<Fo(Gc),f>=<F.(Gc),f>=<F.G,cf>.$$
Hence we have $cA^*\subseteq wap(A)$.\\
iv) The proof is  similar as (iii).\\\\
 Let $G$ be an infinite group. Since $M(G)$ is not
Arens regular, by Theorem 10, the unit element of   $M(G)$ has not
$LW^*W$-property or $RW^*W$-property. Similarly conclusions  hold
for the Banach algebra $l^1(G)$.
  \\\\
{\it{\bf Theorem 11.}}  Let $A$ and $B$ be Banach algebras and let $h$
be a bounded homomorphism
 from $A$ onto $B$. If $A$ has
$LW^*W$-property (resp. $RW^*W$-property), then $B$ has
$LW^*W$-property (resp.
$RW^*W$-property).\\\\
Proof.  Since $h$ is continuous, the second adjoint of $h$,
$h^{**}$ is $weak^*-weak$
 continuous from $A^{**}$ into $B^{**}$.
  Hence we conclude that $h^{**}(F.G)=h^{**}(F).h^{**}(G)$ for
 all $F,G\in A^{**}$.\\
 Let $(g_{\alpha})_{\alpha}\subseteq B^*$ be such that $b.g_{\alpha} \stackrel{w^*} {\rightarrow} 0$ where $b\in
 B$. We define $h^*(g_{\alpha})=f_{\alpha}\in A^*$ and get $h(a)=b$
 for some $a\in A$. Then for every $x\in A$, we have
 $$<a.f_{\alpha},x>=<f_{\alpha},x.a>=<h^*(g_{\alpha}),x.a>=<g_{\alpha},h(x.a)>=<g_{\alpha},h(x).h(a)>$$
 $$<g_{\alpha},h(x).b>=<b.g_{\alpha},h(x)>\rightarrow 0.$$
 Thus $a.f_{\alpha} \stackrel{w^*} {\rightarrow} 0$. Since $a$ has
 $LW^*W$-property, $a.f_{\alpha} \stackrel{w} {\rightarrow} 0$.
 Now, let $G\in B^{**}$. Since $h$ is surjective by [12,
 3.1.22] there is $F\in A^{**}$ such that $h^{**}(F)=G$. Then we
 have the following assertions
 $$<G,b.g_{\alpha}>=<G.b,g_{\alpha}>=<h^{**}(F)h(a), g_{\alpha}>=<h^{**}(F.a), g_{\alpha}>$$
 $$=<(F.a),h^*(g_{\alpha})>=<(F.a), f_{\alpha}>=<F,a. f_{\alpha}>\rightarrow 0$$
Thus we conclude that $b.g_{\alpha} \stackrel{w^*} {\rightarrow}
0$, therefore $B$ has $LW^*W$-property.\\\\
{\it {\bf Theorem 12.}} For a Banach algebra $A$ the
following
assertions hold\\
 i) If $A^*A=A^*$ and $A$ has $RW^*W$-property,
then $A$ is Arens regular.\\
 ii) If $AA^*=A^*$ and $A$ has
$LW^*W$-property, then
$A$ is Arens regular.\\
\\
Proof. i) Let $(F_{\alpha})\subseteq A^{**}$ be a net such that
$F_{\alpha} \stackrel{w^*} {\rightarrow} F$. Then, we have
$(F_{\alpha}.f ).a \stackrel{w^*} {\rightarrow} (F.f).a$ for every
$a\in A$ and $f\in A^*$, and so  $F_{\alpha}.(f .a) \stackrel{w^*}
{\rightarrow} F.(f.a)$. Since $A$ has $LW^*W$-property, we
conclude that $F_{\alpha}.(f .a) \stackrel{w} {\rightarrow}
F.(f.a)$. Therefore for all $G\in A^{**}$, $f\in A^*$ and $a\in A$
we have the following relations
$$<G.F_{\alpha},f.a>=<G,F_{\alpha}(f.a)>=<G,(F_{\alpha}.f).a>\rightarrow
<G,(F.f).a>$$$$=<G.F,f.a>.$$
 Since $A^*A=A^*$, the result follows. \\
 ii) The Proof  is
  similar to (i).\\\\
By Theorem 11 and Theorem 12(ii), if $A$, $B$ are Banach algebras
and $h$ is a bounded
 homomorphism from $A$ onto $B$ whenever $B$  factors on the right, then if $A$ has
 $LW^*W$
 -property, $B$ is Arens regular.
Also from Theorem 12(i) and [11, 2.6] we conclude that when $A$ is
$WSC$ and $AA^*=A^*$ [or $A$ is unital] and $A$ has
$RW^*W$-property then $A$ is Arens regular. \\\\\\
{\it{\bf Lemma 13.}}
For a Banach algebra $A$ the following statements
 hold. \\
i) For all $f\in A^*$ and $F\in A^{**}$ there is a net
$(x_{\alpha})_{\alpha}\subseteq A$
such that $f.x_{\alpha}\stackrel{w^*}{\rightarrow}foF$.\\
ii) $A$ is Arens regular if and only if
$f.x_{\alpha}\stackrel{w}{\rightarrow}foF$ where $F\in A^{**}$,
$f\in A^*$
and $(x_{\alpha})_{\alpha}\subseteq A$ such that $x_{\alpha} \stackrel{w^*} {\rightarrow} F$.
\\\\
Proof. ~i) Since  $w^*-closureA=A^{**}$ , for all $F\in A^{**}$,
there is a net $(x_{\alpha})_{\alpha}\subseteq A$ such that
$x_{\alpha} \stackrel{w^*} {\rightarrow} F$. Now, let $a\in A$.
Then, we have
$$<foF,a>=<F,aof>=\lim_{\alpha}<x_{\alpha},aof>=\lim_{\alpha}<aof,x_{\alpha}>$$$$=\lim_{\alpha}<f,x_{\alpha}.a>
=\lim_{\alpha}<f.x_{\alpha},a>.$$
We conclude that $f.x_{\alpha}\stackrel{w^*} {\rightarrow} foF$.\\
ii) Let $f.x_{\alpha}\stackrel{w} {\rightarrow} foF$. Then, for all $G\in A^{**}$ , we have \\
$$<F.G,f>=<F,G.f>=\lim_{\alpha}<G.f,x_{\alpha}>=\lim_{\alpha}<G,f.x_{\alpha}>=<G,foF>$$$$=<FoG,f>.$$
 Hence $F.G=FoG$ for all $F,G\in A^{**}$.\\
Conversely, let $F.G=FoG$ for all $F,G\in A^{**}$ and let $x_{\alpha}\stackrel{w^*} {\rightarrow} F$.
Then, for all $f\in A^*$ and $G\in A^{**}$, we have the following relations
$$ <G,f.x_{\alpha}>=<G.f,x_{\alpha}>\rightarrow<F,G.f>=<F.G,f>=<FoG,f>$$$$=<G,foF>,$$
which implies  that $f.x_{\alpha}\stackrel{w} {\rightarrow} foF$.\\\\
{\it{\bf Theorem 14.}} For a Banach algebra $A$ the following
assertions
hold.\\
 i)~Suppose~~
$\lim_{\alpha}\lim_{\beta}<f_{\beta},a_{\alpha}>=\lim_{\beta}\lim_{\alpha}
<f_{\beta},a_{\alpha}>$~~ for every\\
$(a_{\alpha})_{\alpha}\subseteq A$ and
$(f_{\beta})_{\beta}\subseteq A^{*}$. Then,
 if $A^*A=A^*$ or $AA^*=A^*$ then  $A$ is Arens regular. Also  if $A^{**}=cA^{**}$
  or  $A^{**}=A^{**}c$ for some $c\in A$, then $A$ is
  reflexive.\\\\
ii)  Let $A$ be Arens regular and let
$(a_{\alpha})_{\alpha}\subseteq A$ and
$(f_{\beta})_{\beta}\subseteq A^{*}$ such that  have taken limits, respectively, in the $weak^*$ and weak topology in $A^{**}$ and $A^*$. Then for all $F\in
  A^{**}$, we have
$$\lim_{\alpha}\lim_{\beta}<F,f_{\beta}.a_{\alpha}>=\lim_{\beta}\lim_{\alpha}<F,f_{\beta}.a_{\alpha}>.$$
iii) If for some $a\in A$, we have
$$\lim_{\alpha}\lim_{\beta}<f_{\beta}a,a_{\alpha}>=\lim_{\beta}\lim_{\alpha}<f_{\beta}a,a_{\alpha}>,$$
then $a$ has $RW^*W$-property. Also if
$$\lim_{\alpha}\lim_{\beta}<af_{\beta},a_{\alpha}>=\lim_{\beta}\lim_{\alpha}<af_{\beta},a_{\alpha}>,$$
then $a$ has $LW^*W$-property.\\\\
Proof. i) If we show that $A$ has  both $LW^*W$-property and
$RW^*W$-property then
 by Theorem 10 and Theorem 12 , we are done.\\
Suppose that $(f_{\beta})_{\beta}\subseteq A^{*}$ and $a\in A$
such that $a.f_{\beta}\stackrel{w^*} {\rightarrow}0$. Let $F\in
A^{**}$. Since~ $w^*-closure A=A^{**}$, there is
$(a_{\alpha})_{\alpha}\subseteq A$ such that
$a_{\alpha}\stackrel{w^*} {\rightarrow}F$. Consequently, we have
the following equalities
$$\lim_{\beta}<F,a.f_{\beta}>=\lim_{\beta}<F.a,f_{\beta}>=\lim_{\beta}\lim_{\alpha}
<a_{\alpha}.a,f_{\beta}>=\lim_{\beta}\lim_{\alpha}<f_{\beta},a_{\alpha}.a>$$$$
=\lim_{\alpha}\lim_{\beta}<f_{\beta},a_{\alpha}.a>
=\lim_{\alpha}\lim_{\beta}<a.f_{\beta},a_{\alpha}>=0.$$
Therefore, we obtain $a.f_{\beta}\stackrel{w} {\rightarrow}0$. Similarly $A$ has $RW^*W$-property.\\
ii) Suppose that $(a_{\alpha})_{\alpha}\subseteq A$ such that
 $a_{\alpha}\stackrel{w^*} {\rightarrow}G$ for some $G\in A^{**}$. Let $(f_{\beta})_{\beta}\subseteq A^{*}$
 such that $f_{\beta}\stackrel{w} {\rightarrow}f$ where $f\in A^*$. Then by Lemma 13,
  for fixed $\beta$, we have $f_{\beta}.a_{\alpha}\stackrel{w} {\rightarrow}f_{\beta}oG$.
 So, for all $F\in A^{**}$ we have the following relations
$$\lim_{\beta}\lim_{\alpha}<F,f_{\beta}.a_{\alpha}>=\lim_{\beta}<F,f_{\beta}.oG>=\lim_{\beta}<GoF,f_{\beta}>=<GoF,f>$$
$$=<G.F,f>=<G,F.f>=\lim_{\alpha}<F.f,a_{\alpha}>=\lim_{\alpha}<F,f.a_{\alpha}>$$
$$=\lim_{\alpha}<a_{\alpha}.F,f>=\lim_{\alpha}\lim_{\beta}<a_{\alpha}.F,f_{\beta}>=\lim_{\alpha}
\lim_{\beta}<F,f_{\beta}.a_{\alpha}>.$$\\
iii) Suppose that $(f_{\beta})_{\beta}\subseteq A^*$ and
$f_{\beta}a\stackrel{w^*}{\rightarrow}0$. Let $F\in A^{**}$ and
$(a_{\alpha})_{\alpha}\subseteq A$ such that
$a_{\alpha}\stackrel{w^*}{\rightarrow}F$. Then, we have
$$\lim_{\beta}<F,f_{\beta}a>=\lim_{\beta}<aF,f_{\beta}>=\lim_{\beta}\lim_{\alpha}<f_{\beta},aa_{\alpha}>$$
$$=\lim_{\alpha}\lim_{\beta}<f_{\beta},aa_{\alpha}>=\lim_{\alpha}\lim_{\beta}<f_{\beta}a,a_{\alpha}>=0$$
Consequently, we have $f_{\beta}a\stackrel{w}{\rightarrow}0$. The prove of the next assertion is the same.\\\\
{\it{\bf Definition 15.}} Let $A$ be Banach algebra. We say that $A^*$ strong left (right) factors, if for all $(f_{\alpha})_{\alpha}\subseteq A^*$, there are $(a_{\alpha})_{\alpha}\subseteq A$ and $f\in A^*$ such that $f_{\alpha}=f.a_{\alpha}~(f_{\alpha}=a_{\alpha}.f)$ where $(a_{\alpha})_{\alpha}\subseteq A$ has limit in the $weak^*$ topology in $A^{**}$. If $A^*$ strong left and right factors, then we say that $A^*$ strong factors.\\
For a Banach algebra $A$ with a BAI, it is clears that if $A^*$ is strong left (right) factors, then $A^*$ factors on the left (right).\\\\
{\it{\bf Theorem 16.}} Let $AA^*\subseteq wapA$. If $A^*$ strong factors on the left (right), then $A$ has $LW^*W$-property ($RW^*W$-property).\\\\
Proof. We prove that $A$ has $LW^*W$-property and the proof of
$RW^*W$-property is similar. Suppose that $a\in A$ and $(f_{\alpha})_{\alpha}\subseteq A^*$ such that $a.f_{\alpha} \stackrel{w} {\rightarrow} 0$. Since $A^*$ strong factors, there are $(a_{\alpha})_{\alpha}\subseteq A$ and $f\in A^*$ such that $f_{\alpha}=f.a_{\alpha}~(f_{\alpha}=a_{\alpha}.f)$ where $(a_{\alpha})_{\alpha}\subseteq A$ has limit in the $weak^*$ topology in $A^{**}$. Let $F\in A^{**}$ and $(a_{\beta})_{\beta}\subseteq A$ such that ${b_{\beta} \stackrel{w^*} {\rightarrow} F} $. Then
$$\lim_{\alpha}<F,a.f_{\alpha}>=\lim_{\alpha}\lim_{\beta}<af_{\alpha},b_{\beta}>=
\lim_{\alpha}\lim_{\beta}<af.a_{\alpha},b_{\beta}>$$
$$=\lim_{\alpha}\lim_{\beta}<af,a_{\alpha}b_{\beta}>=\lim_{\beta}\lim_{\alpha}<af,a_{\alpha}b_{\beta}>
=\lim_{\beta}\lim_{\alpha}<af.a_{\alpha},b_{\beta}>$$
$$=\lim_{\beta}\lim_{\alpha}<af_{\alpha},b_{\beta}>=0.$$
Then $a\in A$ has $LW^*W$-property and so $A$ has $LW^*W$-property.\\\\
{\it{\bf Theorem 17.}} Suppose that $A$ has a $BAI$ and
$AA^{**}=A^{**}$. If $A^*$ is not right factor then there is $a\in
A$ such that $a$  has not
$RW^*W$-property.\\\\
Proof. Let $(e_{\alpha})_{\alpha}$ be a $BAI$ for $A$. Then for
all $f\in A^*$,
 we have $f.e_{\alpha}\stackrel{w^*} {\rightarrow}f$, and so for all $a\in A$,
  we have $(f.e_{\alpha}).a\stackrel{w^*} {\rightarrow}f.a$.
  Suppose on contrary that $a\in A$ has $RW^*W$-property then
   $(f.e_{\alpha}).a\stackrel{w} {\rightarrow}f.a$.\\
Now, for $F\in A^{**}$ we have
$$<a.F,f.e_{\alpha}>=<F,(f.e_{\alpha}).a>=<F.f,e_{\alpha}.a>\rightarrow<F.f,a>=<a.F,f>,$$
hence $f.e_{\alpha}\stackrel{w} {\rightarrow}f$, since
$AA^{**}=A^{**}$. By [11, 2.11] we conclude that
$A^*A=A^*$ which is a contradiction.\\\\
{\it{\bf Open problems }}\\
1. Suppose that $A$ has $LW^*W$-property and $A$ is $WCC$ and
$WSC$. Is $A$ reflexive?\\
2. For a non abelian Banach algebra $A$, when $LW^*W$-property and
$RW^*W$-property are equivalent?\\
 3. If $A^*$  factors on the
right and $A$ has $RW^*W$-property or if $A^*$  factors on the
left and $A$ has $LW^*W$-property,  is
$A$   Arens regular?\\
4. If $A$ factors on the left and Arens regular, is $A$ has $LW^*W$-property?\\
5. If $A$ is  Arens irregular ( or  extremely Arens irregular),
dose $A$ has some member with $LW^*W$-property or
$RW^*W$-property?
\footnotesize{
Department of Mathematic, Amirkabir University of Technology, Tehran, Iran\\
{\it Email address:} haghnejad@aut.ac.ir\\\\
Department of Mathematic, Amirkabir University of Technology, Tehran, Iran\\
{\it Email address:} riazi@aut.ac.ir}}

\end{document}